\documentclass[11pt]{amsart}
\usepackage{amssymb, amsmath}

\title{On equational $\mathfrak{X}$-domains}
\author{O. Al-Raisi}
\address{O. Al-Raisi: Department of  Mathematics,  College of Science, Sultan Qaboos University, Muscat, Oman}
\email{omartalibmiran@gmail.com}

\author{M. Shahryari}
\address{M. Shahryari: Department of  Mathematics,  College of Science, Sultan Qaboos University, Muscat, Oman}
\email{m.ghalehlar@squ.edu.om}

\newtheorem{example}{Example}[section]
\newtheorem{corollary}{Corollary}[section]
\newtheorem{proposition}{Proposition}[section]

\newtheorem {theorem}{Theorem}[section]
\newtheorem{lemma}{Lemma}[section]

\numberwithin{equation}{section}

\newcommand{\XT}{\mathfrak{X}\mathrm{T}}
\newcommand{\CSX}{\mathrm{CS}\mathfrak{X}}
\newcommand{\CT}{\mathrm{CT}}
\newcommand{\CSA}{\mathrm{CSA}}
\newcommand{\X}{\mathfrak{X}}
\newcommand{\Stab}{\mathrm{Stab}}
\newcommand{\Fix}{\mathrm{Fix}}
\newcommand{\Rad}{\mathrm{Rad}}
\begin{document}

\maketitle
\begin{abstract}
Let $\X$ be a variety of groups containing all abelian groups. We introduce and investigate equational $\X$-domains as a generalization of equational domains in algebraic geometry over groups. We characterize these groups by means of systems of equations associated with the identities of $\X$, and describe the radical and coordinate group of the union of two coordinate axes. We also extend classical characterizations of coordinate groups of irreducible algebraic sets, involving discrimination and universal theories, to the setting of equational $\X$-domains. Finally, we establish two construction results: under suitable malnormality, amalgamated free products are equational $\X$-domains; moreover, a restricted regular wreath product of two nontrivial groups is an equational $\X$-domain if and only if its base group is an equational $\X$-domain.
\end{abstract}
\vspace{1cm}


%1111111111111111111111111111111111111111111111111111111111111111111111111
%%%%%%%%%%%%%%%%%%%%%%%%%%%%%%%%%%%%%%%%%%%%%%%%%%%%%%%%%%%%%%%%%%%%%%
\vspace{1cm}

This paper represents a further contribution to our ongoing series of articles  devoted to the development of the theory of $\XT$- and $\CSX$-groups, aimed at extending various notions of algebraic geometry over groups. Recall that a subgroup $H$ of a given group $G$ is called {\em malnormal} if for every element $x\in G\setminus H$, we have  $H\cap H^x=\{1\}$. A group $G$ is called {\em conjugately separated abelian} ($\mathrm{CSA}$) if all maximal abelian subgroups of $G$ are malnormal.
The class of $\mathrm{CSA}$-groups is broad and plays a significant role in the study of residually free groups, universal theory of non-abelian free groups, limit groups, exponential groups and equational domains in algebraic geometry over groups (see \cite{BMR}, \cite{Champ}, and \cite{MR2}). Another class of groups with very close connections to $\mathrm{CSA}$-groups is the class of $\mathrm{CT}$-groups ({\em commutative transitive groups}). A group is $\mathrm{CT}$ if commutativity is a transitive relation on the set of its non-identity elements. Despite this simple definition, the class of CT groups is also central to the study of residually free groups. Every $\mathrm{CSA}$-group is $\mathrm{CT}$ but the converse is not true.

The second author of the present paper proposed an idea in \cite{Shah} for generalizing the notions of $\mathrm{CSA}$- and $\mathrm{CT}$-groups: Suppose that $\X$ is a variety of  groups. A group $G$  is an   $\XT$-group  if and only if, for any two $\X$-subgroups $A$ and $B$ of $G$, the assumption $A\cap B\neq \{1\}$ implies that $\langle A, B\rangle$ is also an $\X$-subgroup of $G$. Similarly,  a group $G$ is a $\CSX$-group if  all of its maximal $\X$-subgroups  are malnormal. The  general theory of $\XT$- and $\CSX$-groups is developed in \cite{Omar-Shah2} and the reader will find the relevant definitions and results there. Given a variety $\X$, a group $G$ and an element $a\in G$, we define the $\X$-centralizer of $a$ in $G$ as the set
$$
C^{\X}_G(a)=\{x\in G:\ \langle a,x \rangle\in \X \}.
$$
Then, as is proved in \cite{Omar-Shah2}, the group $G$ belongs to the class $\XT$ if and only if, for all non-identity elements $a\in G$, either $C^{\X}_G(a)$ is empty or it is an $\X$-subgroup of $G$ . Further, if $\CSX\subseteq \XT$, then the group $G$ belongs to the class $\CSX$ if and only if for every non-identity $a\in G$, the subset  $C^{\X}_G(a)$ is empty or it is a malnormal $\X$-subgroup of $G$.

There is yet another class of groups which play a crucial role in algebraic geometry over groups; the class of equational domains. A group $G$ is called an equational domain (or a domain for short) if for every non-trivial normal subgroup $K$ in $G$ we have $C_G(K)=\{ 1\}$. This is equivalent to saying that for every positive integer $n$ the union of every two algebraic subsets of $G^n$ is an algebraic set. Note that we only consider the case of Diophantine algebraic geometry over $G$ where equations belong to $G[x_1, \ldots, x_n]=G\ast \mathbb{F}_n$, with $\mathbb{F}_n$ being the free group freely generated by the basis $\{ x_1, \ldots, x_n\}$. For details on equational domains, the reader may consult \cite{BMR} and \cite{DMR}. A compact review is also given in \cite{Omar-Shah}. There are close connections between the classes  $\CSA$ (and in general $\CSX$) and equational domains. For example, every non-abelian $\CSA$-group (and every  non-$\X$ group that is a $\CSX$-group) is a domain; consequently, $\CSA$-groups (and more generally, $\CSX$-groups) play an important role in algebraic geometry over groups. This led us to extend the definition of equational domain to a more general concept which we call  {\em equational $\X$-domain} or $\X$-domain in short; the definition is as follows: Consider a variety $\X$  containing all abelian groups. For a group $G$ and a subgroup $K\leq G$, define the $\X$-centralizer of $K$ to be the intersection of all $\X$-centralizers of elements of $K$:
$$
C^{\X}_G(K)=\bigcap_{a\in K}C^{\X}_G(a).
$$
Note that this is not necessarily a subgroup. We say that the group $G$ is an equational $\X$-domain if and only if, for every non-trivial normal subgroup $K\unlhd G$, we have $C^{\X}_G(K)=\{ 1\}$. Note that this is equivalent to each of the following statements:

\begin{itemize}
\item For any elements  $a, b\in G$, if  $\langle a, c\rangle\in \X$ for all $c\in \langle b^G\rangle$, then  $a=1$ or $b=1$.
\item For every non-identity element $b\in G$, we have $C^{\X}_G(\langle b^G\rangle)=\{1\}$.
\end{itemize}

It is quite easy to verify that every $\X$-domain is an ordinary domain, as every abelian group is assumed to belong to $\X$. Also, it can be easily seen that for every variety $\X$ which contains all abelian groups and has the property $\CSX\subseteq \XT$, every $\CSX$-group which is not in $\X$ is an $\X$-domain. Consequently, the classes of $\CSX$-groups, equational domains, and equational $\X$-domains are closely related.

In the first section of this article, we study basic properties of $\X$-domains and their connections to $\CSX$-groups and ordinary domains. The second section concerns the role of $\X$-domains in algebraic geometry over groups and, in particular, it answers this question: What is the necessary and sufficient condition in terms of unions of algebraic sets for a group $G$ to be an $\X$-domain? This section contains general forms of Theorems D2 and D3 in \cite{BMR}, as well. The last section contains two  group-theoretical results which, together, establish the abundance of $\X$-domains. We show, in this last section, that the amalgamated free product of two groups $A$ and $B$ with a common subgroup $C$ that is proper and malnormal in both factors, is an $\X$-domain, even when neither $A$ nor $B$ is assumed to be an $\X$-domain. This generalizes Theorem A2 and Corollary 2 of \cite{BMR} (and of course it modifies these results by removing unnecessary assumptions). It is also proved that a restricted wreath product $A\wr B$ with an equational $\X$-domain group $B$ as the base group, is again an $\X$-domain. As mentioned earlier, the results of this last section demonstrate that the class of $\X$-domains is wide enough to deserve a deep study.

For clarity, we introduce the following notation. Given a subset $X$ of a group $G$, the subgroup generated by $X$ and the normal subgroup generated by $X$ (i.e. the normal closure of $X$ in $G$) will be denoted by $\langle X\rangle$ and $\langle X^G\rangle$, respectively;  a conjugate $a^x$ (or $A^x$) is $x^{-1}ax$ (similarly, $x^{-1}Ax$). A commutator $[x, y]$ is $x^{-1}y^{-1}xy$. The notation $\overline{x}$ stands for the tuple $(x_1, \ldots, x_n)$, and similarly, we use $\overline{a}$ or $\overline{b}$ to denote tuples of elements. The free group of rank $n$ will be denoted by $\mathbb{F}_n$.  We will denote a restricted regular wreath product $A\ltimes B^{(A)}$ by $A\wr B$. {\bf Throughout this paper}, $\X$ will be a variety which {\bf contains all abelian groups}.

\section{Basic properties of $\X$-domains}
As is mentioned in the introduction, every $\X$-domain is an ordinary domain as well. This is because every abelian group belongs to $\X$ and consequently, $C_G(a)\subseteq  C^{\X}_G(a)$ for every $a\in G$. The converse is not true of course as the following example shows.

\begin{example}
Every non-abelian simple group $G$ is a domain; in particular, the alternating group $A_5$ is a domain. Let $\X$ be the variety of groups defined by the law $[x^{60},y]\approx 1$. Note that this variety contains all abelian groups and that $A_5$ is not an $\X$-domain since
$$C^{\X}_{A_5}(A_5)=A_5.$$
\end{example}

It is also easy to see that if $\CSX\subseteq \XT$, then every $\CSX$-group which does not belong to $\X$ is an $\X$-domain. However, for every proper variety $\X$, we can find an $\X$-domain which is not $\CSX$.

\begin{proposition}
If the variety $\X$ is proper, then there is an $\X$-domain which is not $\CSX$.
\end{proposition}

\begin{proof}
Note that every finite group embeds in some alternating group $A_n$ and this means that no proper variety $\X$ contains all finite alternating groups; otherwise $\mathbb{F}_2$ which is residually finite must belong to $\X$. Hence there is an $n\geq 5$ such that $G=A_n\notin \X$. A theorem of Guralnick and Kantor (see \cite{Gur}), and even a very special case of it proved long ago by Sophie Piccard (see \cite{Piccard}), shows that every non-identity element $a\in G$ belongs to a generating pair, i.e. there is another element $b\in G$ such that $G=\langle a, b\rangle$. In order to show that $G$ is an $\X$-domain, we must prove $C^{\X}_G(G)=\{1\}$. Let $a\in C^{\X}_G(G)$ and $b\in G$. Then we have $\langle a, b\rangle \in \X$. But, if $a\neq 1$, then we may choose $b$ such that $G=\langle a, b\rangle$, and consequently $G\in \X$, which is a contradiction. Hence $a=1$, proving that $G$ is an $\X$-domain. However, note that $G$ is not a $\CSX$-group, as we already proved that every finite $\CSX$-group belongs to $\X$ (see \cite{Omar-Shah2}).
\end{proof}

Even if we add the extra assumption of being an $\XT$-group, an $\X$-domain may still not be $\CSX$, as can be verified for the variety $\X$ of all abelian groups and $G=A_5$. Indeed, the alternating group $A_5$ is easily seen to be both a domain and a $\CT$-group; for every non-identity element $a\in G$, the centralizer $C_G(a)$ is isomorphic to one of the abelian groups $C_3, C_5, V_4\in \X$. However, it is not a $\CSA$-group as all finite $\CSA$-groups are known to be abelian.

The following is a list of some basic facts about $\X$-domains. We don't provide a proof here because the proofs of these basic facts mirror the proofs of the corresponding facts for ordinary equational domains.

\begin{proposition}
 Let $G$ be a non-trivial $\X$-domain. Each of the following is true:
 \begin{enumerate}
     \item The group $G$ does not belong to the class $\X$.
     \item Every normal subgroup of $G$ which belongs to the class $\X$ is trivial; as such, $G$ is not solvable, and hence, not nilpotent.
     \item The group $G$ is directly indecomposable.
 \end{enumerate}
\end{proposition}

It is easy to see that a non-trivial finite group is an equational domain if and only if it is monolithic with a non-abelian monolith. Likewise, for the case of finite $\X$-domains, we have the following:

\begin{theorem}
Non-trivial finite $\X$-domains are exactly the finite monolithic groups such that their monolith  does not belong to $\X$.
\end{theorem}

\begin{proof}
Let $G$ be a non-trivial finite $\X$-domain. Then $G$ is an ordinary equational domain. Let $K$ and $L$ be two distinct minimal normal subgroups of $G$. Then $[K, L]=1$ as their intersection is trivial, and this implies that $L\subseteq C_G(K)$, which is a contradiction. Hence, $G$ is monolithic. Suppose $M$ is the monolith of $G$. If $M\in \X$, then $C^{\X}_G(M)=M$, another contradiction. This proves one direction. 

Conversely, suppose that $G$ is a finite monolithic group with monolith $M\notin\X$. Since $M$ is a minimal normal subgroup of $G$ and $\X$ contains all abelian groups,
$$
M=T_1\times\cdots\times T_k,
$$
where all factors  are isomorphic to a non-abelian finite simple group $T$ which does not belong to $\X$. First, note that  $C_G(M)=1$. Indeed, if $C_G(M)\ne1$, then by uniqueness of the minimal normal subgroup, $M\subseteq  C_G(M)$, contradicting $Z(M)=1$. Fix $1\neq x\in G$. Conjugation by $x$ permutes the simple direct factors of $M$. Since $C_G(M)=\{1\}$, there is an orbit $I$ of this permutation such that $x$ induces a nonidentity automorphism $\alpha$ of $N=\prod_{i\in I}T_i$. Identify $N$ with $\mathrm{Inn}(N)$, the group of inner automorphisms of $N$, and put $$
Q=\langle\mathrm{Inn}(N),\alpha\rangle\subseteq \mathrm{Aut}(N).
$$
The subgroup $\mathrm{Inn}(N)\cong N$ is a minimal normal subgroup of $Q$: Any normal subgroup contained in it is a product of some of its
simple factors, and $\alpha$ permutes those factors transitively. Also $C_Q(\mathrm{Inn}(N))=\{1\}$, because $N$ is centerless and
$$
C_{\mathrm{Aut}(N)}(\mathrm{Inn}(N))=\{1\}.
$$
Consequently, every nontrivial normal subgroup $D$ of $Q$ contains $\mathrm{Inn}(N)$. Indeed, otherwise $D\cap\mathrm{Inn}(N)=\{1\}$,
so $D$ centralizes $\mathrm{Inn}(N)$, a contradiction. Thus every proper quotient of $Q$ is cyclic, since
$$
\frac{Q}{\mathrm{Inn}(N)}=\langle\alpha\mathrm{Inn}(N)\rangle.
$$

By the theorem of Burness, Guralnick, and Harper (see \cite{BGH}) on generating pairs in finite groups whose proper quotients are cyclic, there exists
$\beta \in Q$ such that $Q=\langle\alpha,\beta\rangle$. Since $Q=\mathrm{Inn}(N)\langle\alpha\rangle$, we may write $\beta=m\alpha^r$ for some $m\in N$ and integer $r$. Hence $Q=\langle\alpha,m\rangle$. The conjugation action of $\langle x,m\rangle$ on $N$ therefore maps onto $Q$. Since $Q$ contains $N\cong T^{|I|}\notin\X$, it follows that $\langle x,m\rangle\notin\X$. 

Thus no non-identity $x\in G$ belongs to $C^{\X}_G(M)$, so $C^{\X}_G(M)=1$. Every nontrivial normal subgroup of $G$ contains $M$; therefore $G$ is an $\X$-domain.
\end{proof}

Note that the above  argument proves a slightly stronger statement. If $M$ is a non-abelian minimal normal subgroup of an arbitrary finite group $G$
with $M\notin\X$, then $C^{\X}_G(M)=C_G(M)$. Indeed, the proof shows that $x\notin C_G(M)$ implies $x\notin C^{\X}_G(M)$, while the reverse inclusion follows because $\X$ contains all abelian groups. It is not difficult to show that every  monolithic group with non-abelian monolith is an equational domain. However, the case of $\X$-domains seems to be hard: Suppose $G$ is a monolithic group with the monolith not in $\X$. Can we conclude that $G$ is an $\X$-domain? Unfortunately, at the moment, we have neither a proof, nor a counterexample.

As we saw,  $\X$-domains within the class of $\XT$-groups are not necessarily $\CSX$. However, it will soon be apparent that such $\X$- domains are easier to work with. Let's call an element $x\in G$ a {\em zero divisor} relative to $\X$, if it belongs to $C^{\X}_G(K)$ for some non-trivial normal subgroup $K$. So,  being $\X$-domain is equivalent to non-existence of non-trivial zero-divisors relative to $\X$.

\begin{proposition}
Let $G$ be an $\XT$-group. A non-identity element $x$ of $G$ is a zero-divisor relative to $\X$ if and only if there is a non-identity element $y$ of $G$ such that $$\langle x,y^g\rangle\in \X$$ for all $g\in G$.
\end{proposition}

\begin{proof}
Suppose there is a pair of non-identity elements $x,y$ of $G$ such that $\langle x,y^g\rangle\in \X$ for all $g\in G$. Using the $\XT$-property of $G$, it is easy to see that
$$
x\in C^{\X}(\langle y^G\rangle).
$$
Indeed, an arbitrary element of $\langle y^G\rangle$ is of the form
$$
g_1^{-1}y^{k _1}g_1\cdots g_m^{-1}y^{k_m}g_m
$$
where $k_1,\ldots, k_m$ are integers and $g_1,\ldots,g_m$ are in $G$. Now, from $\langle x,g_i^{-1}yg_i\rangle\in \X$, we infer that the subgroup $\langle x,g_i^{-1}y^{k_i}g_i\rangle$ of $\langle x,g_i^{-1}yg_i\rangle$ belongs to $\X$ for each $i$; accordingly, the $\XT$-property implies that
$$
\langle x,g_1^{-1}y^{k_1}g_1,\ldots,g_m^{-1}y^{k_m}g_m \rangle\in \X.
$$
Therefore,
$$
\langle x, g_1^{-1}y^{k_1}g_1\cdots g_m^{-1}y^{k_m}g_m \rangle \in \X
$$
which proves our claim that $x\in C^{\X}(\langle y^G\rangle)$. Conversely, if $x\in G$ is a zero-divisor relative to $\X$ then there is a non-identity element $y$ of $G$ such that $x\in C^{\X}(\langle y^G\rangle)$. In particular, this implies that $\langle x,y^g\rangle\in \X$ for all $g\in G$.
\end{proof}

Consequently, we have the following characterization of $\XT$-groups that are $\X$-domains.

\begin{corollary}
A necessary and sufficient condition for an $\XT$-group $G$ to be an $\X$-domain is that for every $x,y\in G$, if $\langle x,y^g\rangle\in \X$ for all $g\in G$, then $x=1$ or $y=1$.
\end{corollary}

We close this section with the following:

\begin{proposition}
Consider an $\XT$-group $G$ which contains a non-trivial normal subgroup $N$ that is an $\X$-domain. Then $G$ is an  $\X$-domain if and only if
$C^{\X}_G(N) = \{1\}$.
\end{proposition}

\begin{proof}
Assume that $C^{\X}_G(N) = \{1\}$. Let $U$ be a non-trivial normal subgroup of $G$; we claim that $U\cap N\neq \{1\}$. Indeed, if $U\cap N=\{1\}$ then, the normality of $U$ in $G$ implies that, $[U,N]\subseteq U\cap N=\{1\}$; this means that $\{1\}\neq U\subseteq C_G(N)\subseteq C^{\X}_G(N)$ which is a contradiction. Therefore, every non-trivial normal subgroup of $G$ must intersect $N$ non-trivially.
Set $B=U\cap N$ and note that $B$ is a normal subgroup of $G$ which is contained in $N$. Since $N$ is an $\X$-domain, we must have $C_N^{\X}(B)=C_G^{\X}(B)\cap N=\{1\}$. But  $C_G^{\X}(B)$ is a normal subgroup of $G$ and hence, $C^{\X}_G(B)=\{1\}$. Now, $B\subseteq U$ so $C_G^{\X}(U)\subseteq C_G^{\X}(B)$ which forces the desired conclusion that $C_G^{\X}(U)=\{1\}$.
\end{proof}

\section{Algebraic geometry over $\X$-domains}
To review the basic concepts of algebraic geometry over groups, we ask the reader to consult \cite{BMR}.  Here, we fix our notations. Let $G$ be a group. A $G$-group is a group $H$ which contains $G$ as a subgroup. Suppose $X=\{ x_1, \ldots, x_n\}$ is the set of variables. Then $G[X]$ is the free product of the group $G$ and the free group $\mathbb{F}[X]$. For an element $p(\overline{x})\in G[X]$, the corresponding equation is denoted by $p(\overline{x})\approx 1$. For a system of equations $S\subseteq G[X]$, the corresponding algebraic set in $G^n$ and its radical  will be denoted by $V_G(S)$ and $\Rad_G(S)$, respectively. Also the notation $\Rad(Y)$ will be used for the radical of a subset $Y\subseteq G^n$. The notation $\Gamma_G(S)$ (respectively, $\Gamma(Y)$) will be used for the coordinate group of the system $S$ (the subset $Y$). Note that these coordinate groups are finitely generated $G$-groups.

Recall that a group $G$ is an equational domain if and only if, for every positive integer $n$, the union of two algebraic subsets of $G^n$ is an algebraic set: More precisely, as is proved in \cite{BMR}, when $G$ is an equational domain then, for any systems $S, T\subseteq G[X]$, we have
$$
V_G(S)\cup V_G(T)=V_G(\{ [p(\overline{x}), q(\overline{x})^g]:\ p(\overline{x})\in S, q(\overline{x})\in T, g\in G\}).
$$
Furthermore, $G$ is an equational domain if and only if
$$
A=\{ (a, b)\in G^2:\ a=1\ \mathrm{or}\ b=1\}=V_G(x\approx 1)\cup V_G(y\approx 1)
$$
is an algebraic set.

The discussion in the previous paragraph motivates us to ask a question: Is there any similar description of equational $\X$-domains? In this section  we answer this question. As before, $\X$ is an arbitrary variety which contains all abelian groups. Define a subset of $G^2$ as follows:
\begin{eqnarray*}
Z_{\X}(G)&=&\{ (a, b)\in G^2:\ b\in C^{\X}_G(\langle a^G\rangle)\}\\
         &=&\{ (a, b)\in G^2:\  \langle b, c\rangle \in \X\ \quad \forall c\in \langle a^G\rangle\}.
\end{eqnarray*}
It is easy to see that $V_G(x\approx 1)\cup V_G(y\approx 1)\subseteq Z_{\X}(G)$. The reverse inclusion holds exactly when $G$ is an $\X$-domain:

\begin{lemma}\label{Z-1}
The set $Z_{\X}(G)$ is an algebraic set. The group $G$ is an $\X$-domain if and only if $V_G(x\approx 1)\cup V_G(y\approx 1)=Z_{\X}(G)$.
\end{lemma}

\begin{proof}
Note that every element $c\in \langle a^G\rangle$ has the form $c=v(a^{g_1}, \ldots, a^{g_k})$ for some elements $g_1, \ldots, g_k\in G$ and a group word $v(x_1, \ldots, x_k)$. Hence
$$
\langle b, c\rangle =\{ u(b, v(a^{g_1}, \ldots, a^{g_k})):\ u(x, y)\in \mathbb{F}_2\}.
$$
Hence, $\langle b, c\rangle \in \X$ if and only if, for all $m\geq 1$, all $m$-variable identities $w\approx 1$ of the variety $\X$, and all $u_1, \ldots, u_m\in \mathbb{F}_2$, we have
$$
w(u_1(b, v(a^{g_1}, \ldots, a^{g_k})), \ldots, u_m(b, v(a^{g_1}, \ldots, a^{g_k})))=1.
$$
So, consider the set $\Sigma_{\X}$ which consists of the equations of the form
$$
w(u_1(y, v(x^{g_1}, \ldots, x^{g_k})), \ldots, u_m(y, v(x^{g_1}, \ldots, x^{g_k})))\approx 1,
$$
where $m, k\geq 1$, $g_1, \ldots, g_k\in G$, $u_1, \ldots, u_m\in \mathbb{F}_2$, $v\in \mathbb{F}_k$, and $w\approx 1$ is an $m$-variable law of $\X$. Then we have $Z_{\X}(G)=V_G(\Sigma_{\X})$, hence it is an algebraic set.
\end{proof}

Consequently, we have

\begin{corollary}\label{Sigma-1}
The group $G$ is an equational $\X$-domain if and only if
$$
V_G(\Sigma_{\X})=V_G(x\approx 1)\cup V_G(y\approx 1).
$$
\end{corollary}

In order to generalize this result for the union of two arbitrary algebraic sets, we need to introduce more special systems of equations which rely on the identities of the variety $\X$. Suppose $p(\overline{x}), q(\overline{x})\in G[X]$ and consider the set
$$
Z_{\X}(p, q)=\{ \overline{a}\in G^n:\ q(\overline{a})\in C^{\X}_G(\langle p(\overline{a})^G\rangle)\}.
$$

\begin{lemma}\label{Z-2}
If the group $G$ is an $\X$-domain, then
$$
Z_{\X}(p, q)=V_G(p(\overline{x})\approx 1)\cup V_G(q(\overline{x})\approx 1).
$$
\end{lemma}

\begin{proof}
As $\X$ contains all abelian groups, it is obvious that
$$
V_G(p(\overline{x})\approx 1)\cup V_G(q(\overline{x})\approx 1)\subseteq Z_{\X}(p, q).
$$
For the reverse inclusion, let $\overline{a}\in Z_{\X}(p, q)$. Then the pair $(p(\overline{a}), q(\overline{a}))$ belongs to $Z_{\X}(G)$ and  hence, by Lemma \ref{Z-1}, we have
$$
(p(\overline{a}), q(\overline{a}))\in V_G(x\approx 1)\cup V_G(y\approx 1).
$$
This means that $\overline{a}\in V_G(p(\overline{x})\approx 1)\cup V_G(q(\overline{x})\approx 1)$ which completes the proof.
\end{proof}

Now, substitute $p$ and $q$ in every element of $\Sigma_{\X}$ for $x$ and $y$, respectively, and denote the new set by $\Sigma_{\X}(p, q)$. Then we have

\begin{corollary}\label{Sigma-2}
The group $G$ is an equational $\X$-domain if and only if for all $p(\overline{x}), q(\overline{x})\in G[X]$,
$$
V_G(\Sigma_{\X}(p, q))=V_G(p(\overline{x})\approx 1)\cup V_G(q(\overline{x})\approx 1).
$$
\end{corollary}

Therefore, for any two  systems $S, T\subseteq G[X]$, consider the set
$$
\Sigma_{\X}(S, T)=\bigcup_{p\in S, q\in T}\Sigma_{\X}(p, q).
$$
Note that we always have
$$
V_G(S)\cup V_G(T)=\bigcap_{p\in S, q\in T}(V_G(p(\overline{x})\approx 1)\cup V_G(q(\overline{x})\approx 1)),
$$
and consequently, we have:

\begin{theorem}\label{Alg-set}
The group $G$ is an $\X$-domain if and only if for any systems $S, T\subseteq G[X]$, the equality
$$
V_G(\Sigma_{\X}(S, T))=V_G(S)\cup V_G(T)
$$
holds.
\end{theorem}

In the case when $\X$ is the variety of all abelian groups, the set $\Sigma_{\X}$ is equivalent to
$$
C(x, y)=\{ [x, y^g]:\ g\in G\},
$$
and hence, $G$ is an equational domain if and only if
$$
V_G(C(x,y))=V_G(x\approx 1)\cup V_G(y\approx 1).
$$
This means that, $G$ is an equational domain if and only if it satisfies the (possibly infinite) sentence
$$
\forall x\ \forall y \left(\bigwedge_{g\in G}([x, y^g]\approx 1)\rightarrow (x\approx 1 \vee y\approx 1)\right).
$$
For an arbitrary variety which contains all abelian groups, if $G$ is an $\X$-domain, then it is an ordinary domain as well, so we have
$$
V_G(C(x,y))=V_G(\Sigma_{\X}),
$$
which is equivalent to saying that $G$ satisfies the (possibly infinite) sentence
$$
\forall x\ \forall y \left(\bigwedge_{g\in G}([x, y^g]\approx 1)\leftrightarrow \bigwedge_{w\in \Sigma_{\X}}(w(x, y)\approx 1)\right).
$$
It also says that if $G$ is an $\X$-domain, then $\Rad_G(C(x, y))=\Rad_G(\Sigma_{\X})$. However, note that the converse implication is not true: If  $\Rad_G(C(x, y))=\Rad_G(\Sigma_{\X})$, then one may not be able to conclude $G$ is an $\X$-domain. For example, suppose $G$ is a non-trivial abelian group. Obviously, $G$ is not an $\X$-domain, but we have
$$
V_G(C(x, y))=G^2=V_G(\Sigma_{\X}).
$$

Our aim is to determine the radical of $\Sigma_{\X}$, explicitly, and to find its coordinate group, in the case when $G$ is an $\X$-domain. Returning to the set
$$
A=V_G(x\approx 1)\cup V_G(y\approx 1),
$$
first we determine the coordinate group of $A$. Note that
$$
\Rad(G)=\{ p(z)\in G[z]:\ p(a)=1\ \forall a\in G\}.
$$
Hence, let's denote the coordinate group of $G$ by
$$
P_G=\frac{G[z]}{\Rad(G)}.
$$
We denote the coset containing $p(z)$ by the same notation $p(z)$.  Let $\varepsilon: P_G\to G$ be defined as $\varepsilon(z)=1$ and $\varepsilon(g)=g$, for all $g\in G$. Now, define the fiber product
$$
P_G\times_{\varepsilon}P_G=\{ (p(z), q(z))\in P_G\times P_G:\ \varepsilon(p(z))=\varepsilon(q(z))\}.
$$
Put $J_G=\Rad(A)$.

\begin{lemma}\label{Coor-A}
We have $\Gamma(A)=P_G\times_{\varepsilon}P_G$.
\end{lemma}

\begin{proof}
Recall that
$$
\Gamma(A)=\frac{G[x, y]}{J_G}.
$$
Define a homomorphism $f:G[x, y]\to P_G\times_{\varepsilon}P_G$ by extending the assignments
$$
f(x)=(z, 1), \ f(y)=(1, z), \ f(g)=(g, g)\ \forall g\in G.
$$
Let $q(x, y)\in \ker(f)$. Then $(q(z, 1), q(1, z))=(1, 1)$ and hence,
$$
q(z, 1), q(1, z)\in \Rad(G).
$$
Therefore, for all $a\in G$ we have $q(a, 1)=q(1, a)=1$, and this shows that $q(x, y)\in J_G$. A similar argument shows that $J_G\subseteq \ker(f)$. To show that $f$ is surjective, let $(p, q)\in P_G\times_{\varepsilon}P_G$. Suppose $g=\varepsilon(p)=\varepsilon(q)$. Then $p(1)=q(1)=g$, and if we consider the element $q(x, y)=p(x)g^{-1}q(y)$, then

\begin{eqnarray*}
f(q(x, y))&=&(p(z)g^{-1}q(1), p(1)g^{-1}q(z))\\
          &=& (p(z), q(z)).
\end{eqnarray*}
This shows that $G[x, y]/J_G\cong P_G\times_{\varepsilon}P_G$.
\end{proof}

Let $K$ be the kernel of the homomorphism $\varepsilon$. Note that if $g\in G\cap K$, then $g=\varepsilon(g)=1$, and consequently, $G\cap K=\{1\}$. For every $p(z)\in P_G$, we have $p(z)=gg^{-1}p(z)$, where $g=p(1)$, and we have $\varepsilon(g^{-1}p(z))=1$. This shows that $P_G$ is the semidirect product of $G$ and $K$. Hence, we proved:

\begin{corollary}\label{Coor-A-2}
We have $\Gamma(A)\cong G\ltimes (K\times K)$, where $G$ acts on $K\times K$ diagonally as
$$
(p(z), q(z))^g=(g^{-1}p(z)g, g^{-1}q(z)g).
$$
\end{corollary}
Hence, we conclude the following:

\begin{theorem}\label{Coor-Sigma}
Let $G$  be an $\X$-domain. Then $\Rad_G(\Sigma_{\X})=J_G$ and
$$
\Gamma_G(\Sigma_{\X})\cong G\ltimes(K\times K).
$$
\end{theorem}

The next goal is to describe the radical $J_G$ explicitly. Suppose
$$
I_G(z)=\Rad(G)=\{ p(z)\in G[z]:\ p(a)=1\ \forall a\in G\}.
$$
By replacing $z$ with $x$ and $y$, we obtain two subgroups of $G[x, y]$ which we denote them by $I_G(x)$ and $I_G(y)$, respectively.

\begin{lemma}\label{Ncl-1}
The radical $J_G$ is the normal closure of the set
$$
I_G(x)\cup I_G(y)\cup \{ [x, y^g]:\ g\in G\}
$$
in the group $G[x, y]$.
\end{lemma}

\begin{proof}
Suppose $N$ is the normal closure of the set given above and note that $N\subseteq J_G$. To prove the reverse inclusion, consider the group $Q=G[x, y]/N$. Note that the group $G[x]/I_G(x)$ is an isomorphic copy of $P_G$, where $z$ is replaced by $x$, so we denote this quotient group by $P_G(x)$. Similarly, $G[y]/I_G(y)=P_G(y)$. Now we have the presentation
\begin{eqnarray*}
Q&=&\langle G[x], G[y]:\ I_G(x), I_G(y), [x, y^g]\ (g\in G)\rangle\\
 &=&\langle P_G(x), P_G(y):\ [x^h, y^g]\ (g, h\in G)\rangle.
\end{eqnarray*}
Note that $P_G(x)$ and $P_G(y)$ share $G$ as an amalgamating subgroup and therefore
$$
Q=\langle P_G(x)\ast_{_G}P_G(y):\ [x^h, y^g]\ (g, h\in G)\rangle.
$$
Recall that we proved $P_G=G\ltimes K$ and hence, $P_G(x)=G\ltimes K_x$, $P_G(y)=G\ltimes K_y$, where $K_x$ and $K_y$ are the copies of $K$ after replacing $z$ with $x$ and $y$, respectively. Hence
$$
Q=\langle (G\ltimes K_x)\ast_{_G}(G\ltimes K_y):\ [x^h, y^g]\ (g, h\in G)\rangle.
$$
One verifies easily that $K_x=\langle x^g:\ g\in G\rangle$ and $K_y=\langle y^g:\ g\in G\rangle$. Hence, in $Q$, every conjugate $K_x^{\ h}$ commutes with every conjugate $K_y^{\ g}$. This means that $Q=G\ltimes(K_x\times K_y)\cong G\ltimes (K\times K)$. Consequently,
$$
\frac{G[x, y]}{N}\cong \frac{G[x, y]}{J_G}.
$$
Now, recall the epimorphism $f:G[x, y]\to P_G\times_{\varepsilon}P_G$ from the proof of Lemma \ref{Coor-A}. As $N\subseteq J_G$, there is an induced map
$$
\overline{f}:\frac{G[x, y]}{N}\to P_G\times_{\varepsilon}P_G,
$$
such that $\ker(\overline{f})=J_G/N$. But, $G[x, y]/N=Q\cong G\ltimes (K_x\times K_y)$ and we have
$$
\overline{f}(xN)=(z, 1), \ \overline{f}(yN)=(1, z), \ \overline{f}(gN)=g\  (\forall g\in G).
$$
This shows that $\overline{f}$ is injective, and hence, $N=J_G$.
\end{proof}

\begin{corollary}\label{Ncl-2}
A group $G$ is an $\X$-domain if and only if $\Rad_G(\Sigma_{\X})$ is the normal closure of the set
$$
I_G(x)\cup I_G(y)\cup \{ [x, y^g]:\ g\in G\}
$$
in the group $G[x, y]$.
\end{corollary}

The article \cite{BMR} (Theorems D2 and D3) contains a characterization of coordinate groups of irreducible algebraic sets over equational domains. Recall that a $G$-group is a group $H$ which contains a distinguished copy of $G$. A $G$-group $H$ is said to be separated by $G$ (or in other words, $H$ is residually $G$), if for  every non-identity element $h\in H$ there is a homomorphism $\varphi:H\to G$ such that the restriction $\varphi_{|_G}$ is the identity map and $\varphi(h)\neq 1$. Similarly, we say that $H$ is discriminated by $G$ (or fully residually $G$), if for every finite number of non-identity elements $h_1, \ldots, h_m$, there is a homomorphism $\varphi:H\to G$ such  that the restriction $\varphi_{|_G}$ is the identity map and none of the elements $\varphi(h_1), \ldots, \varphi(h_m)$ equals the identity of $G$. A $G$-group $H$ is called an equational $G$-domain if for any $a, b\in H$, the assumption $[a, b^G]=1$ implies $a=1$ or $b=1$. Finally, a normal subgroup $P\unlhd H$ is called $G$-prime subgroup if $H/P$ is an equational $G$-domain. In \cite{BMR}, the following result is proved (we re-state it in the case of Diophantine geometry, only):

\begin{theorem}\label{D2-3}
Let $G$ be  an equationally Noetherian domain and $Y\subseteq G^n$ be a non-empty algebraic set. Then the following are equivalent:
\begin{enumerate}
\item The algebraic set $Y$  is irreducible.
\item The subgroup $\Rad(Y)$ is a $G$-prime subgroup of $G[X]$.
\item The coordinate group $\Gamma(Y)$ is a $G$-equationally Noetherian $G$-domain.
\item The coordinate group $\Gamma(Y)$ is discriminated by $G$.
\item The groups $G$ and $\Gamma(Y)$ have the same universal theories in the first order language of groups with parameters in $G$.
\end{enumerate}
\end{theorem}

Our aim is to formulate an analogue of the above theorem for the case of $\X$-domains. If $H$ is a $G$-group, we say that it is a $(G, \X)$-domain if and only if, for every non-identity element $a\in H$, we have $C^{\X}_H(\langle a^G\rangle)=\{1\}$.

\begin{lemma}
Let $G$ be an $\X$-domain and $H$ be a $G$-group which is separated by $G$. Then the following are equivalent:
\begin{enumerate}
\item $H$ is discriminated by $G$.
\item $H$ is a $(G, \X)$-domain.
\item $H$ is an $\X$-domain.
\end{enumerate}
\end{lemma}

\begin{proof}
Note that every $(G, \X)$-domain is a $G$-domain, i.e. for all $a, b\in H$, we have
$$
[a, b^G]=1\rightarrow a=1\ \mathrm{or}\ b=1.
$$
So, we can apply Theorem C1 of \cite{BMR} to prove the implication $(2)\rightarrow (1)$: The group $H$ is $G$-domain, $G$ is a domain, and $H$ is separated by $G$, hence $H$ is discriminated by $G$. To prove $(1)\rightarrow (2)$, suppose $H$ is discriminated by $G$ but it is  not $(G, \X)$-domain. Hence, there are non-identity elements $a, b\in H$ such that $b\in C^{\X}_H(\langle a^G\rangle)$. So, for all $u\in \langle a^G\rangle$ we have $\langle u, b\rangle \in \X$. By $(1)$, there is a homomorphism $\varphi:H\to G$ such that both $\varphi(a)$ and $\varphi(b)$ are non-identity and $\varphi$ fixes the elements of $G$. Note that
$$
\varphi(\langle a^G\rangle)=\langle \varphi(a)^G\rangle.
$$
Let $v$ be an element from the right hand side. Then $v=\varphi(u)$ for some $u\in \langle a^G\rangle$. As we have $\langle u, b\rangle \in \X$, by applying $\varphi$, we get
$$
\langle v, \varphi(b)\rangle =\varphi(\langle u, b\rangle)\in \X.
$$
Hence $\varphi(b)\in C^{\X}_G(\langle \varphi(a)^G\rangle)$. But, $G$ is an $\X$-domain, so $\varphi(b)=1$, a contradiction.

Now, we prove the implication $(2)\rightarrow (3)$. Let $K$ be a non-trivial normal subgroup of $H$ and $b\in C^{\X}_H(K)$. Let $a\in K$ be a non-identity element, so $\langle a^G\rangle \subseteq K$ and consequently, $b\in C^{\X}_H(\langle a^G\rangle)$. But, $H$ is a $(G, \X)$-domain, so $b=1$, proving that $H$ is an $\X$-domain. (Note that the separation assumption is not used here.)

Finally, to prove $(3)\rightarrow (2)$, let $H$ be an $\X$-domain, but not $(G, \X)$-domain. So, there are non-identity elements $a, b\in H$ such that $b\in C^{\X}_H(\langle a^G\rangle)$. Note that $C^{\X}_H(\langle a^H\rangle)=1$ as $H$ is an $\X$-domain, hence $b$ does not belong to $C^{\X}_H(\langle a^H\rangle)$. Therefore, there is a $c\in \langle a^H\rangle$ such that $\langle b, c\rangle$ does not belong to $\X$. This means that there is a law in the variety $\X$, say $w(x_1, \ldots, x_m)\approx 1$, such that it fails in $\langle b, c\rangle$. Hence, there are elements $u_1, \ldots, u_m\in \langle b, c\rangle$ such that
$$
w(u_1, \ldots, u_m)\neq 1.
$$
Now, since $H$ is separated by $G$, there is a homomorphism $\varphi:H\to G$ such that
$$
\varphi(w(u_1, \ldots, u_m))\neq 1.
$$
We claim that $\varphi(b)\neq 1$ and $\varphi(c)\neq 1$. If $\varphi(c)=1$, then all the elements $\varphi(u_1), \ldots, \varphi(u_m)$ belong to the cyclic group $\langle \varphi(b)\rangle$ and hence
$$
\varphi(w(u_1, \ldots, u_m))=1
$$
as all abelian groups belong to $\X$. Thus $\varphi(c)\neq 1$ and a similar argument shows that $\varphi(b)\neq 1$. Recall that $c\in \langle a^H\rangle$, so $\varphi(a)\neq 1$, as well. Hence, we have
$$
\varphi(b)\in C^{\X}_G(\langle \varphi(a)^G\rangle)
$$
which violates the fact that $G$ is an $\X$-domain.
\end{proof}

In what follows, by a prime subgroup, we mean a normal subgroup $P\unlhd G$ such that $G/P$ is an equational domain. Similarly, $P$ is called $\X$-prime if $G/P$ is an $\X$-domain. Now, we are ready to formulate the $\X$-version of Theorem \ref{D2-3}.

\begin{theorem}\label{X-version}
Let $G$ be an equationally Noetherian $\X$-domain. Let $Y\subseteq G^n$ be a non-empty algebraic set. Then the following are equivalent:
\begin{enumerate}
\item The algebraic set $Y$ is irreducible.
\item The subgroup $\Rad(Y)$ is an $\X$-prime subgroup of $G[X]$.
\item The coordinate group $\Gamma(Y)$ is a $G$-equationally Noetherian $\X$-domain.
\item The coordinate group $\Gamma(Y)$ is discriminated by $G$.
\item The groups $G$ and $\Gamma(Y)$ have the same universal theories in the first order language of groups with parameters in $G$.
\end{enumerate}
\end{theorem}

\begin{proof}
Note that the coordinate group $\Gamma(Y)$ is separated by $G$, so we can apply the previous lemma for the $G$-group $H=\Gamma(Y)$. As every $\X$-domain is an ordinary domain, the equivalence of $(1)$ and $(4)$ comes from Theorem \ref{D2-3}. The equivalence of $(3)$ and $(4)$ is a result of the previous lemma and Theorem \ref{D2-3}. Also, as $\Gamma(Y)$ is separated by $G$, it embeds into a direct power of $G$ and consequently, it is $G$-equationally Noetherian, so we conclude the equivalence of $(2)$ and $(3)$. Finally, the equivalence of $(1)$ and $(5)$ comes from Theorem \ref{D2-3}.
\end{proof}

Recall from \cite{BMR} that if $H$ is a finitely generated $G$-group with the same universal  theory as $G$ (with parameters in $G$), then $H$ is the coordinate group of some irreducible algebraic set $Y\subseteq G^n$ (for some $n\geq 1$), and hence it is an $\X$-domain.

\section{A free construction and wreath product}
There are several methods to construct new equational domains: In \cite{BMR} it is shown that the amalgamated free product of two equational domains with a malnormal amalgam is a domain. A similar result for the case of $\CSA$-groups is proved in \cite{Gilden}. We extended the same result for $\CSX$-groups in \cite{Omar-Shah4}. It is also possible to construct equational domains using wreath products: In \cite{BMR} it is proved that a restricted wreath product $A\wr B$ is an equational domain if $A$ and $B$ are domains and $A$ is torsion-free. This result is modified by us in \cite{Omar-Shah3} by removing the assumptions on $A$: The wreath product $A\wr B$ is an equational domain if and only if $B$ is an equational domain. In this section, we prove analogous results for the case of $\X$-domains. First we discuss the amalgamated free product $G=A\ast_{_C}B$. The details of our argument are very similar to that in \cite{Omar-Shah4} and consequently, we don't repeat those technical details here. Since our argument is rely on Bass-Serre theory, we review the Bass-Serre tree of $G=A\ast_{_C}B$:

Define an undirected graph $T$ as follows: $V_T$, the set of vertices, is the collection of all left cosets of $A$ and $B$ in $G$ i.e.,
$$
V_T=\{ gA, gB:\ g\in G\}.
$$
The set of edges, $E_T$, is the collection of all pairs $\{ gA, gB\}$ with $g$ in $G$, in other words
$$
E_T=\{ \{ gA, gB\}:\ g\in G\}.
$$
According to the main theorem of Bass-Serre theory, this graph is in fact a tree. The group $G$ acts on it by  left multiplication and the vertex and edge stabilizers are
$$
\Stab_G(gA)=gAg^{-1},\ \Stab_G(gB)=gBg^{-1},\ \Stab_G(\{ gA, gB\})=gCg^{-1}.
$$
The reader may use \cite{Dicks} or \cite{Serre} for details on this tree.

\begin{theorem}\label{Amal-dom}
Let $\X$ be a proper variety containing all abelian groups and suppose $A$ and $B$ are groups without any involutions. Let $C$ be a proper malnormal subgroup of $A$ and $B$. Then $G=A\ast_{_C}B$ is an equational $\X$-domain.
\end{theorem}

\begin{proof}
Let $T$ be the Bass-Serre tree of $G$. As is shown in \cite{Omar-Shah4}, every non-identity element of $G$ fixes at most one edge of $T$ and consequently, for every $1\neq x\in G$, the set $\Fix_T(x)$ (the set of all fixed points of $x$) contains at most two vertices. Let $K$ be a non-trivial normal subgroup of $G$. \\

{\bf Claim i}. $K$ contains a hyperbolic element.\\

To prove the claim, let $1\neq x\in K$ be an arbitrary element and suppose $\Fix_T(x)$ is non-empty. Let $a\in A\setminus C$ and $b\in B\setminus C$. Then the element $t=ab\in G$ is cyclically reduced of length two and hence $t$ is hyperbolic, so it acts as a translation along its axis. Note that, for any positive integer $n$, $\Fix_T(t^nxt^{-n})=t^n\Fix_T(x)$. For sufficiently large $n$, we have
$$
\Fix_T(x)\cap t^n\Fix_T(x)=\emptyset
$$
and hence by Lemma 1.2 of \cite{Omar-Shah4}, the element $y=x(t^nxt^{-n})\in K$ is hyperbolic. \\

{\bf Claim ii}. Any two hyperbolic elements with different axes generate a non-abelian free group after taking sufficiently large powers.\\

To see why, suppose $x$ and $y$ are hyperbolic elements of $G$ with the axes $L_x$ and $L_y$, respectively. Let $L_x\neq L_y$. There are three cases:\\

{\bf i}. Suppose $L_x\cap L_y$ contains an infinite ray. Let $m$ and $n$ be the translation lengths of $x$ and $y$, respectively and consider the elements $x^n$ and $y^m$. They are both translations of the same translation length $mn$. If needed we can replace $y$ by $y^{-1}$ to make sure that $x^n$ and $y^m$ translate the common ray in the same direction. Hence, for any vertex $w$ on the common ray, we have $x^nw=y^mw$, and so $(y^{-m}x^n)w=w$. As this happens for infinitely many vertices on the common ray, and as every non-identity element of $G$ fixes at most two vertices, we conclude that $x^n=y^m$. Consequently, the axes of $x^n$ and $y^m$ are equal. But, note that the axis of $x^n$ is  $L_x$ and the axis of $y^m$ is $L_y$. This contradiction shows that the intersection $L_x\cap L_y$ is a finite set. \\

{\bf ii}. Suppose $L_x\cap L_y$ is non-empty but finite. In this case, we use a similar argument as in \cite{Omar-Shah4} to prove that for sufficiently large $m$ and $n$, $\langle x^m, y^n\rangle \cong \mathbb{F}_2$.\\

{\bf iii}. Let $L_x\cap L_y=\emptyset$. Again using a similar argument as in \cite{Omar-Shah4}, we see that $\langle x, y\rangle\cong \mathbb{F}_2$.\\

{\bf Claim iii}. The normal subgroup $K$ contains two hyperbolic elements with different axes.\\

To prove this claim, suppose $x\in K$ is hyperbolic. Let $L$ be the axis of $x$ and choose a vertex $v=gA$ in $L$. Note that the degree of $v$ is the index $[A:C]\geq 3$. So, there is at least one edge of $T$, incident to $v$ and not in $L$. We know that the stabilizer $\mathrm{Stab}_G(v)$ acts transitively on the set of edges incident to $v$. Hence, there is an element $g\in \mathrm{Stab}_G(v)$ which sends an incident edge of $v$ in  $L$ to another incident edge of $v$ which is not in $L$. This means that $gL\neq L$. Now, let $y=gxg^{-1}$ which is an element of $K$. This new element is also hyperbolic and it has an axis different from the axis of $x$.\\

{\bf Claim iv}. For every $1\neq g\in G$, there is a $h\in K$ such that $\mathbb{F}_2\subseteq \langle g, h\rangle$.\\

We prove this claim by considering two cases:\\

{\bf i}. Suppose $g$ is hyperbolic and $L_g$ is its axis. By Claim iii above, there are $x, y\in K$, both hyperbolic with different axes $L_x$ and $L_y$, respectively. Hence, at least one of these axes (say $L_x$) is not equal to $L_g$. Consequently, $\mathbb{F}_2\subseteq \langle g, x\rangle$, by Claim ii. \\

{\bf ii}. Suppose $g$ has a fixed vertex. Choose $x\in K$ hyperbolic with the axis $L$. Suppose $gL=L$. Then $g$ must act as a translation on the line $L$ as $G$ does not contain any involution. But, $g$ fixes a vertex, so it must fix the whole $L$. This implies that $g=1$. Therefore, $gL\neq L$. Hence by Claim ii, we have $\mathbb{F}_2\subseteq \langle gxg^{-1}, x\rangle\subseteq \langle g, x\rangle$. \\

{\bf Claim v}. We have $C^{\X}_G(K)=\{1\}$.\\

This is true, because if we take $1\neq g\in C^{\X}_G(K)$, then by Claim iv, there is an element $h\in K$ such that $\mathbb{F}_2\subseteq \langle g, h\rangle$. But, at the same time, we have $\langle g, h\rangle \in \X$, and consequently $\mathbb{F}_2\in \X$, which is a contradiction.
\end{proof}

Note that in the above theorem, one does not need to assume $A$ and $B$ are $\X$-domains; we only need two arbitrary groups without any involutions which share a proper malnormal subgroup. Also, the theorem shows that $G=A\ast_{_C}B$ is an $\X$-domain for every proper variety $\X$ which contains all abelian groups.

The next result is another construction which employs wreath product.
Let $A$ and $B$ be non-trivial groups and consider the restricted regular wreath product $A\wr B=A\ltimes B^{(A)}$. Here, $B^{(A)}$ denotes the group of all functions $f:A\to B$ with finite supports and  $A$ acts on $B^{(A)}$ as follows: For any $a, u\in A$ and each $f\in B^{(A)}$, we have $f^a(u)=f(au)$. An arbitrary element of $G=A\wr B$ is of the form $af$, where  $a\in A$, $f\in B^{(A)}$, and the product of two elements $a_1f_1$ and $a_2f_2$ of $G$ is given by
$$
(a_1f_1)(a_2f_2)=(a_1a_2)(f_1^{a_2}f_2).
$$
The inverse of $af$ is given by $a^{-1}(f^{-1})^{a^{-1}}$, where $f^{-1}$ is the pointwise inverse of $f$.

\begin{theorem}\label{wreath}
Let $A$ and $B$ be non-trivial groups and $\X$ be a variety of groups containing all abelian groups. Then the restricted  wreath product $G=A\wr B=A\ltimes B^{(A)}$ is an $\X$-domain if and only if $B$ is an $\X$-domain.
\end{theorem}

\begin{proof}
For each $u\in A$, we shall write $\pi_u:B^{(A)}\to B$ for the homomorphism $f\mapsto f(u)$; for each $b\in B$ and $u\in A$, we shall write $f_{u,b}$ for the map $A\to B$ which sends $u$ to $b$ and sends every $u'\in A\setminus\{u\}$ to $1\in B$.
Suppose the wreath product is an $\X$-domain but the base group $B$ is not. Then $B$ contains a non-trivial normal subgroup $K$ and a non-identity element $x$ such that $\langle x,y\rangle\in \X$ for all $y\in K$. Let
$$
K_1=K^{(A)}=\{ h\in B^{(A)}:h(u)\in K \ (\forall u\in A)\}.
$$
Note that $K_1\neq\{1\}$ because $K\neq \{1\}$; for example, if $y$ is a non-identity element of $K$ then the map $f_{1, y}$ belongs to $K_1$. Furthermore, $K_1$ is a normal subgroup of $G$ because, for every $g=af$ in $G$ and every $h\in K_1$, the map
$$
(af)^{-1}h(af)=f^{-1}h^af:A\to B
$$
sends every $u \in A$ to $f^{-1}(u)h(au)f(u)$ which belongs to $K$ as $h(au)\in K$, $f(u)\in B$, and $K$ is normal in $B$.

Consider the map $c=f_{1,x}$ and let $m\in K_1$, and put $H=\langle c,m\rangle\leq B^{(A)}$. Suppose that $\{u_1,\ldots, u_r\}$ is the support of $m$ with the point $1$ removed (so every $u_j \neq 1$, and possibly $r=0$), and set $S=\{1,u_1,\ldots,u_r\}$. Since elements of $H$ are products of the elements $c^{\pm1}, m^{\pm1}$, the support of every element of $H$ is a subset of $S$.  Now
$$
\pi_1(H)=\langle c(1),m(1)\rangle=\langle x,m(1)\rangle\in \X,
$$
and, for all $j\in \{1,\ldots,r\}$,
$$
\pi_{u_j}(H)=\langle 1,m(u_j)\rangle=\langle m(u_j)\rangle\in \X
$$
because $\X$ contains all abelian groups. Therefore,
$$
\pi_1(H)\times \pi_{u_1}(H)\times\cdots \times\pi_{u_r}(H)\in \X.
$$
The map
$$
\phi:H\to \pi_1(H)\times \pi_{u_1}(H)\times\cdots \times\pi_{u_r}(H)
$$
given by $\phi(h)=(h(1),h(u_1),\ldots,h(u_r))$ is a monomorphism. It follows that $H$ belongs to the variety $\X$; that is, for each $m\in K_1$, we have $\langle c,m\rangle\in \X$. This shows that $C^{\X}_G(K_1)$ contains a non-identity element of $G$ which is $c$. Thus, $G$ is not an $\X$-domain.

For the converse, suppose that $G$ is not an $\X$-domain so that $G$ has a non-trivial normal subgroup $K$ such that $C^{\X}_G(K)$ contains a non-identity element $g=af$, where $a\in A$ and $f\in B^{(A)}$. Our goal is to extract a normal subgroup $K_0$ of $B$ such that $C_B^{\X}(K_0)\neq \{1\}$.
We begin by observing that, for the normal subgroup $K$ of $G$ as described above, the intersection $K\cap B^{(A)}$ is not trivial. To see this, suppose  $b$ is a non-trivial element of $B$ and $y=af$ be an element of $K$. If $a=1$ then we are done, so we may assume otherwise. Because both $K$ and $B^{A}$ are normal in $G$,
$$
[f_{1,b},y]=(f_{1,b})^{-1}y^{-1}f_{1,b}y\in K\cap B^{(A)}
$$
and
$$
((f_{1,b})^{-1}y^{-1}f_{1,b}y)(1)=((f_{1,b})^{-1}f^{-1}f_{1,b}^af)(1)=b^{-1}\neq 1
$$
shows that $[f_{1,b},y]\neq 1$. Let $m$ be a non-identity element of $K\cap B^{(A)}$, say at $u_0\in A$, we have $m(u_0)\neq 1$. If the center $Z(B)$ of $B$ is not trivial then $B$ is not an $\X$-domain and we are done. Thus, we may assume that $Z(B)=\{1\}$, so that there is an element $b_0\in B$ with $[m(u_0),b_0]\neq 1$. It is easy to see that
$$
[m,f_{u_0,b_0}]=m^{-1}((f_{u_0,b_0})^{-1}mf_{u_0,b_0})\in K\cap B^{(A)},
$$
and $[m,f_{u_0,b_0}]=f_{u_0,[m(u_0),b_0]}\neq 1$. We now have a natural candidate for our normal subgroup $K_0$ of $B$: Put
$$
K_0=\{ b \in B: f_{u_0,b}\in K \};
$$
the argument given above shows that $K_0$ contains non-trivial elements of $B$, and it is easy to see that $K_0$ is a normal subgroup of $B$. We claim that $f_{u,b}\in K$ for all $u\in A$ and all $b\in K_0$. Indeed, $K$ is normal in $G$, $f_{u_0,b}\in K$, and
$$
f_{u,b}=(f_{u_0,b})^{u_0u^{-1}}.
$$
We now finish the proof arguing by cases.\\

{\bf i}. Suppose $a=1$ so that $g=f$. There is a $u\in A$ such that $f(u)\neq1$; claim that $f(u)\in C^{\X}_B(K_0)$. To see this, recall from the previous paragraph that, for each $b\in K_0$, we have $f_{u,b}\in K$. This means that $\langle f,f_{u,b}\rangle\in \X$ for every $b\in K_0$. But then, the homomorphic image
$$
\pi_u(\langle f,f_{u,b} \rangle)=\langle f(u),b\rangle
$$
of $\langle f,f_{u,b} \rangle)$ belongs to $\X$ for all $b\in K_0$ which is the desired conclusion.\\

{\bf ii}. Suppose that $a\neq 1$; we claim that $K_0\subseteq C^{\X}_B(K_0)$. To see this, fix $b_0\in K_0\setminus \{1\}$ and let $b\in K_0$. We will prove that $\langle b,b_0\rangle\in \X$. The element $f(1)bf(1)^{-1}$ belongs to $K_0$ since $K_0$ is normal in $B$. Define $h\in B^{(A)}$ by $h(1)=b_0$, $h(a)=f(1)bf(1)^{-1}$, and $h=1$ elsewhere; clearly,
$$
h=f_{1,b_0}f_{a,f(1)bf(1)^{-1}}\in K
$$
because each of $f_{1,b_0}$ and $f_{a,f(1)bf(1)^{-1}}$ belongs to $K$ by the paragraph preceding Case i. This implies that $\langle g,h\rangle\in \X$.
Now,
$$
\langle h,g^{-1}hg\rangle\leq\langle h,g\rangle\in \X
$$
implies that $\langle h,g^{-1}hg\rangle\in \X$. But then
$$
\pi_1(\langle h,g^{-1}hg\rangle)=\langle h(1),(g^{-1}hg)(1)\rangle=\langle b_0,b\rangle\in \X.
$$
This completes the proof.
\end{proof}

As an application, we have the following corollary.\\

\begin{corollary}
Let $\X$ be a proper variety of groups containing all abelian groups. Then $\mathbb{F}_2 \wr \mathbb{F}_2$ is an $\X$-domain which does not belong to $\mathrm{CS}\mathfrak{Y}$ for any proper variety $\mathfrak{Y}$ containing all abelian groups.
\end{corollary}

\begin{proof}
By the previous result, given any proper variety $\X$ which contains all abelian groups, the group $G=\mathbb{F}_2 \wr \mathbb{F}_2$ is an $\X$-domain because,  $\mathbb{F}_2$ is an $\X$-domain by Theorem \ref{Amal-dom}. Assume that $G\in \mathrm{CS}\mathfrak{Y}$ for some variety of groups $\mathfrak{Y}$ containing all abelian groups. Then, as $\mathrm{CS}\mathfrak{Y}$ is a subgroup-closed class of groups (see Proposition 1.4 in \cite{Omar-Shah2}), $\mathbb{F}_2^{(\mathbb{F}_2)}\in \mathrm{CS}\mathfrak{Y}$. Now, if $\mathbb{F}_2^{(\mathbb{F}_2)}\notin\mathfrak{Y}$ then by Theorem 4.1 in \cite{Omar-Shah2}, $\mathbb{F}_2^{(\mathbb{F}_2)}$ would be a domain which is impossible as domains are directly indecomposable. Therefore, $\mathbb{F}_2^{(\mathbb{F}_2)}\in \mathfrak{Y}$. But then $\mathbb{F}_2\in \mathfrak{Y}$ which implies that $\mathfrak{Y}$ is the variety of all groups.
\end{proof}

\end{document}